\documentclass[a4paper,12pt,twoside]{article}

\title{The Mixture of Markov Jump Processes: Monte Carlo Method and the EM Estimation}

 \author{H. Frydman\footnote{ NYU Leonard N. Stern School of Business, Kaufman Management Center 44 West Fourth Street, 8-55 New York, NY 10012, USA. Email address: hfrydman@stern.nyu.edu} \; and\; B.A. Surya\footnote{School of Mathematics and Statistics, Victoria University of Wellington, Gate 6 Kelburn PDE, Wellington 6140, New Zealand. Email address: budhi.surya@vuw.ac.nz }\\ Department of IOMS-Statistics\\
 New York University Stern School of Business\\ and\\
 School of Mathematics and Statistics \\ Victoria University of Wellington}

\date{}

\usepackage{latexsym,epic,eepic,float}
\usepackage{color,colordvi}
\usepackage{shadow}
\usepackage{amsmath,amssymb,amsfonts,stackrel}
\usepackage{amsmath}
\usepackage{enumitem}
\usepackage{epsfig}
\usepackage{lscape}
\usepackage{rotating}
\usepackage{subfigure}
\usepackage[export]{adjustbox}
\usepackage{multicol}
\usepackage{tikz}
\usepackage{graphicx}
\usepackage{bbold}
\usetikzlibrary{automata,positioning}

\newtheorem{theorem}{Theorem}[section]
\newtheorem{theo}[theorem]{Theorem}
\newtheorem{lem}[theorem]{Lemma}

\newtheorem{prop}[theorem]{Proposition}
\newtheorem{Rem}[theorem]{Remark}

\newcommand{\exit}{{\mbox{\, \vspace{3mm}}} \hfill\mbox{$\square$}}

\rm 

\numberwithin{equation}{section}

\date{27 November 2018}

\begin{document}

\maketitle \pagestyle{myheadings} \markboth{H. Frydman and B.A. Surya}{The Mixture of Markov Processes: Monte Carlo and the EM Estimation}

\begin{abstract}
This paper discusses tractable development and statistical estimation of a continuous time stochastic process with a finite state space having non-Markov property. The process is formed by a finite mixture of right-continuous Markov jump processes moving at different speeds on the same finite state space, whereas the speed regimes are assumed to be unobservable. The mixture was first proposed by Frydman \cite{Frydman2005} and recently generalized in Surya (\cite{Surya2018b},\cite{Surya2018}), in which distributional properties and explicit identities of the process are given in its full generality. The contribution of this paper is two fold. First, we present Monte Carlo method for constructing the process and show distributional equivalence between the simulated process and the actual process. Secondly, we perform statistical inference on the distribution parameters of the process. Under complete observation of the sample paths, consistent maximum likelihood estimations are given in explicit form in terms of sufficient statistics of the process. Estimation under incomplete observation is performed using the EM algorithm. The estimation results completely characterize the process in terms of the initial probability of starting the process in any phase of the state space, intensity matrices of the underlying Markov jump processes, and the switching probability matrix of the process. Some numerical examples are given to test the performance of the developed method. The proposed estimation generalizes statistical inferences for the Markov model \cite{Albert}, the mover-stayer model \cite{Frydman1984} and the restricted Markov mixture model \cite{Frydman2005}. 

\medskip

\textbf{MSC2010 Subject Classification:} 60J20, 60J27, 60J28, 62N99


\end{abstract}

\section{The mixture of Markov jump processes}

Throughout the remaining of this paper we denote by $X=\{X^{(\phi)}(t), t\geq0\}$ a Markov mixture process, which is a continuous-time stochastic process defined as a finite mixture of independent Markov jump processes $X^{(m)}=\{X^{(m)}(t): t\geq 0\}$, with $m=1,\dots,M$, whose intensity matrices are given by $\{\mathbf{Q}^{(m)}\}$. We assume that the underlying Markov processes $\{X^{(m)}\}$ have right-continuous sample paths, defined on the same state space $\mathbb{S}=\{1,\dots,p\}$. It is defined by
\begin{equation*}
X=
\begin{cases}
X^{(1)}, & \phi=1\\
\vdots \\
X^{(M)}, & \phi=M
\end{cases}
\end{equation*}
where the variable $\phi$ represents the speed regimes, assumed to be unobservable. This is to say that when the realization of the mixture process $X$ is observed, we do not know from which speed regime $\phi$ the observed process came from. 

More conveniently, we can represent $X$ in terms of the underlying processes $\{X^{(m)}\}$ as follows. Define a Bernoulli indicator variable $\Phi^{(m)}=\mathbb{1}_{\{\phi=m\}}$, which will be used later for the estimation, see (\ref{eq:data}). Notice that $\sum_{m=1}^M \Phi^{(m)}=1.$ Thus,
\begin{equation}\label{eq:mixture}
X(t)=\sum_{m=1}^M \Phi^{(m)} X^{(m)}(t) \quad \textrm{for $t\geq 0$}.
\end{equation}

It is clear that $X$ (\ref{eq:mixture}) represents a finite mixture of Markov processes $X^{(m)}$, and that the random variable $\Phi^{(m)}$ may in general depend on the realization of $X$ as $\Phi^{(m)}=1$ if and only if $X=X^{(m)}$. This implies that the conditional probability $\mathbb{P}\{\Phi^{(m)}=1\vert X(s), 0\leq s \leq t\}$ depends on the past realizations of the process. Empirical evidence of this fact can be found in \cite{Frydman2008}. We refer to \cite{Surya2018b} and \cite{Surya2018} for further distributional properties and explicit identities of the mixture process, in particular in the presence of stochastically closed (absorbing) sets. 

The entry $\{q_{ij}^{(m)}:i,j=1,\dots,p\}$ of matrix $\mathbf{Q}^{(m)}$ satisfies the properties:
\begin{equation}\label{eq:matq}
q_{ii}^{(m)}\leq 0, \; \; q_{ij}^{(m)}\geq 0, \; \; \sum_{j\neq i} q_{ij}^{(m)}=-q_{ii}^{(m)}=q_i^{(m)}, \quad (i,j)\in \mathbb{S}.
\end{equation}

For a given initial state $i_0\in\mathbb{S}$, there is a separate mixing probability
\begin{equation}\label{eq:portion}
s_{i_0}^{(m)}=\mathbb{P}\{\phi=m \vert X_0=i_0\} \quad \textrm{with}\quad \sum_{m=1}^M s_{i_0}^{(m)}=1,
\end{equation}
and $0\leq s_{i_0}^{(m)} \leq 1$. The quantity $s_{i_0}^{(m)}$ has the interpretation as the proportion of population with initial state $i_0$ evolving w.r.t to $X^{(m)}$. In general, $X^{(k)}$ and $X^{(l)}$, $k\neq l$, have different expected length of occupation time of a state $i$, i.e., $1/q_{i}^{(k)}\neq 1/q_{i}^{(l)}$, and have different probability of leaving state $i\in E$ to state $j\in\mathbb{S}$, $j\neq i$, i.e. $q_{ij}^{(k)}/q_{i}^{(k)}\neq q_{ij}^{(l)}/q_{i}^{(l)}$. Note that we have used $q_i^{(m)}$ and $q_{ij}^{(m)}$ to denote the negative of the $i$th diagonal element and the $(i,j)$ entry of $\mathbf{Q}^{(m)}$.

%--------------------------
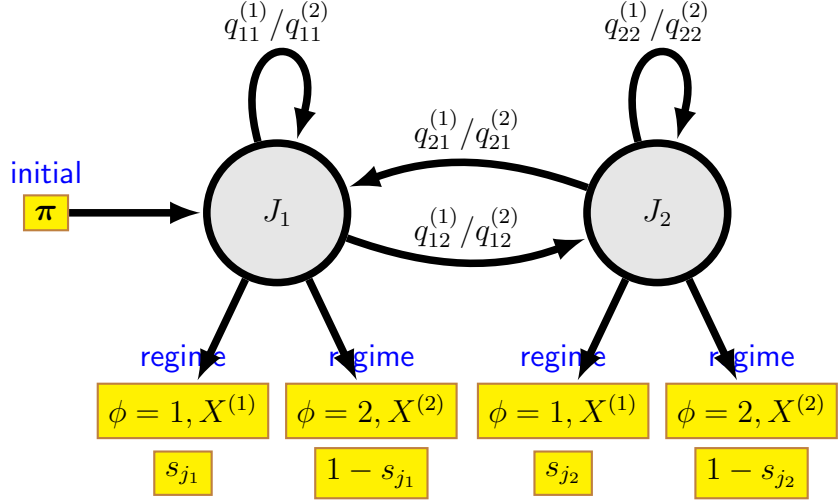
\begin{figure}
\begin{center}
  \begin{tikzpicture}[font=\sffamily]

        % Setup the style for the states
        \tikzset{node style/.style={state,
                                    minimum width=1.85cm,
                                    line width=1mm,
                                    fill=gray!20!white}}

          \tikzset{My Rectangle/.style={rectangle, draw=brown, fill=yellow, thick,
    prefix after command= {\pgfextra{\tikzset{every label/.style={blue}}, label=below}}
    }
}

          \tikzset{My Rectangle2/.style={rectangle,draw=brown,  fill=yellow, thick,
    prefix after command= {\pgfextra{\tikzset{every label/.style={blue}}, label=below}}
    }
}

          \tikzset{My Rectangle3/.style={rectangle, draw=brown, fill=yellow, thick,
    prefix after command= {\pgfextra{\tikzset{every label/.style={blue}}, label=below}}
    }
}

        \node[node style] at (2, 0)     (s1)     {$J_1$};
        \node[node style] at (7, 0)     (s2)     {$J_2$};

      \node [My Rectangle3, label={initial}] at  ([shift={(-5em,0em)}]s1.west) (p0) {$\boldsymbol{\pi}$};

        \node [My Rectangle, label={regime} ] at  ([shift={(3em,-4em)}]s1.south) (g1) {$\phi=2, X^{(2)}$};
          \node [My Rectangle, label={regime} ] at  ([shift={(-3em,-4em)}]s1.south)  (g2) {$\phi=1, X^{(1)}$};

          \node [My Rectangle2] at  ([shift={(3em,-6em)}]s1.south) {$1-s_{j_1}$};
          \node [My Rectangle2] at  ([shift={(-3em,-6em)}]s1.south) {$s_{j_1}$};

         \node [My Rectangle, label={regime} ] at  ([shift={(3em,-4em)}]s2.south) (g3) {$\phi=2, X^{(2)}$};
          \node [My Rectangle, label={regime} ] at  ([shift={(-3em,-4em)}]s2.south)  (g4) {$\phi=1, X^{(1)}$};

           \node [My Rectangle2] at  ([shift={(3em,-6em)}]s2.south) {$1-s_{j_2}$};
          \node [My Rectangle2] at  ([shift={(-3em,-6em)}]s2.south) {$s_{j_2}$};

        % Connect the states with arrows
        \draw[every loop,
              auto=right,
              line width=1mm,
              >=latex,
              draw=orange,
              fill=orange]

            (s1)  edge[bend right=20, auto=left] node {$q_{12}^{(1)}/q_{12}^{(2)}$} (s2)
            (s1)  edge[loop above]                     node {$q_{11}^{(1)}/q_{11}^{(2)}$} (s1)
            (s2)  edge[loop above]                     node {$q_{22}^{(1)}/q_{22}^{(2)}$} (s2)
            (s2)  edge[bend right=20]                node {$q_{21}^{(1)}/q_{21}^{(2)}$}  (s1)

            (s1) edge node {} (g1)
            (s1) edge node {} (g2)

            (s2) edge node {} (g3)
            (s2) edge node {} (g4)

            (p0) edge node {} (s1);

 \end{tikzpicture}
 \caption{State diagram of the Markov mixture process (\ref{eq:mixture}) with $m=2$.}\label{fig:mixture}
\end{center}
\end{figure}

%%-----------------------

Figure \ref{fig:mixture} illustrates the transition of $X$ for the mixture of two Markov jump processes moving from state $J_1$ to $J_2$, and vice versa. When $X$ is observed in state $J_1$, it would stay in the state for an exponential period of time with intensity $q_{j_1}^{(1)}$ or $q_{j_1}^{(2)}$ before moving to $J_2$ with probability $q_{j_1,j_2}^{(1)}/q_{j_1}^{(1)}$ or $q_{j_1,j_2}^{(2)}/q_{j_1}^{(2)}$, depending on whether it is either driven by the underlying Markov process $X^{(1)}$ or $X^{(2)}$.

Markov mixture process is a generalization of mover-stayer model, a mixture of two discrete-time Markov chains proposed in 1955 by Blumen et al \cite{Blumen} to model population heterogeneity in jobs labor market. In the mover-stayer model \cite{Blumen}, the population of workers consists of stayers (workers who always stay in the same job category, $\mathbf{Q}^{(1)}=\mathbf{0}$) and movers (workers who move to other job according to a stationary Markov chain with intensity matrix $\mathbf{Q}^{(2)}$). Estimation of the mover-stayer model was discussed in Frydman \cite{Frydman1984}. Frydman \cite{Frydman2005} generalized the model to a finite mixture of Markov chains moving with different speeds. Frydman and Schuermann \cite{Frydman2008} later on used the result for the mixture of two Markov jump processes moving with intensity matrices $\mathbf{Q}^{(1)}$ and $\mathbf{Q}^{(2)}=\boldsymbol{\Psi}\mathbf{Q}^{(1)}$, i.e., $q_{ij}^{(2)}=\psi_i q_{ij}^{(1)}$, where $\boldsymbol{\Psi}$ is a diagonal matrix, to model the dynamics of firms' credit ratings. Depending on whether $0=\psi_i:=[\boldsymbol{\Psi}]_{i,i}$,  $0<\psi_i<1$, $\psi_i>1$ or $\psi_i=1$, $X^{(2)}$ never moves out of state $i$ (the mover-stayer model), moves out of state $i$ at lower rate, higher rate or at the same rate, subsequently, than that of $X^{(1)}$. If $\psi_i=1$, for all $i\in \mathbb{S}$, $X$ reduces to a simple Markov jump process $X^{(1)}$. However, the mixture model considered in \cite{Frydman2005} and \cite{Frydman2008} is restricted to the case in which the underlying Markov processes $X^{(1)}$ and $X^{(2)}$ have the same probability of leaving a state to another different state. That is $q_{ij}^{(2)}/q_i^{(2)}=q_{ij}^{(1)}/q_i^{(1)}$, $j\neq i$.

The mixture process $X$ has appealing features that, unlike its underlying process $X^{(m)}$, the mixture itself lacks the Markov property; future development of its state depends on its past information and the current time. We refer to \cite{Frydman2008}, \cite{Surya2018b} and \cite{Surya2018} for further distributional properties of the mixture process. 

The transition probability matrix $\mathbf{P}(t)$ of $X$ is given following \cite{Frydman2005} and \cite{Surya2018b} by
\begin{equation}\label{eq:transM}
\mathbf{P}(t)= \sum_{m=1}^M \mathbf{S}^{(m)} e^{\mathbf{Q}^{(m)}t} \quad \mathrm{with} \quad \sum_{m=1}^M \mathbf{S}^{(m)} =\mathbf{I}, 
\end{equation}
for all $t\geq0$, where $\mathbf{I}$ is $(p\times p)-$identity matrix, whereas $\mathbf{S}^{(m)}$ denotes a $(p\times p)-$diagonal matrix, representing the switching probability matrix of $X$, i.e.,
\begin{equation}\label{eq:St}
\mathbf{S}^{(m)} =
 \left(\begin{array}{cc}
 s_1^{(m)} & \mathbf{0} \\
  \mathbf{0} & s_{p}^{(m)} \\
\end{array}\right).
\end{equation}

It is clear from (\ref{eq:transM}) and (\ref{eq:St}) that the distribution of the mixture process $X$ (\ref{eq:mixture}) is uniquely characterized by the variables $\{q_{ij}^{(m)}\}$, $\{q_{i}^{(m)}\}$, respectively representing the off and diagonal elements of $\mathbf{Q}^{(m)}$, the element $\{s_i^{(m)}\}$ of switching probability matrix $\mathbf{S}$, and the probability distribution $\{\pi_i\}$ of starting $X$ in any state $i\in\mathbb{S}$. Furthermore, when we set $\mathbf{Q}^{(m)}=\mathbf{Q}$, all underlying Markov processes move at the same speed $\mathbf{Q}$, $X$ becomes just a $\mathbf{Q}-$Markov process.

\section{Construction of the mixture process}

This section discusses construction of the mixture process $X$ (\ref{eq:mixture}), which can be used to generate the sample paths of the mixture process using Monte Carlo method developed by adapting the approach of Sections 2.1 and 5.1 in Resnick \cite{Resnick}. The simulated sample paths will later be used to solve the inverse problem of estimating the distribution parameters of the process, given its full or incomplete observation of the sample paths, which is the subject of Sections 3 and 4.   

\subsection{Finite mixture of Markov chains}

To start with, let $U_0,U,V,$ and $W$ be independent uniform $[0,1]$ random variables. Introduce a discrete-time Markov chain $Z^{(m)}=\{Z_n^{(m)}\;:n\in\mathbb{N}_+\}$, living on the same finite state space $\mathbb{S}$. It is defined as the corresponding embedded Markov chain for the Markov process $X^{(m)}$. The transition probability matrix of $Z^{(m)}$ is specified by $(p\times p)-$matrix $\boldsymbol{\Pi}^{(m)}$ whose $(i,j)-$element is defined by
\begin{align}\label{eq:MatrM}
\pi_{ij}^{(m)}=
\begin{cases}
q_{ij}^{(m)}/q_i^{(m)}, & j\neq i\\
0, & j=i.
\end{cases}
\end{align}
In the mixture model \cite{Frydman2005}, each embedded Markov chain $Z^{(m)}$ was assumed to have the same transition probability matrix $\boldsymbol{\Pi}$, i.e., $\boldsymbol{\Pi}^{(m)}=\boldsymbol{\Pi}$ for all $m=1,\dots,M$. 

Assume $X$ chooses its initial state $X_0=i_0$ randomly with probability $\boldsymbol{\pi}$:
\begin{align}\label{eq:initialcond}
X_0=\sum_{k=1}^p k \mathbb{1}_{[\sum_{i=1}^{k-1} \pi_i, \; \sum_{i=1}^{k} \pi_i)}(U_0), 
\end{align}
where, we set $\sum_{i=1}^0 \pi_i=0$. Applying similar idea, the speed regime $\phi=m$ of $X$ can be selected randomly, given initial state $X_0=i_0$, at probability $s_{i_0}^{(m)}$ using  
\begin{align}\label{eq:initialregime}
\phi =\sum_{m=1}^M m  \mathbb{1}_{[ \sum_{k=1}^{m-1} s_{i_0}^{(k)}, \; \sum_{k=1}^{m} s_{i_0}^{(k)})}(U).
\end{align}

In the sequel below we denote respectively by $\{V_n\}$ and $\{W_n\}$ $n$ independent copies of the random variables $V$ and $W$ independent of $U_0$ and $U$. The result below gives a Monte Carlo construction of the discrete-time Markov chain $Z_n^{(m)}$.

\begin{lem}\label{lem:mainlemma}
For a given $m=1,\dots,M$, the process $\{Y_n:n\in\mathbb{N}_+\}$ defined by
\begin{equation}\label{eq:chainmixture2}
\begin{split}
Y_{n+1}=& \sum_{k=1}^p k  \mathbb{1}_{ [\sum_{j=1}^{k-1} [\boldsymbol{\Pi}^{(m)}]_{Y_n,j}, \; \sum_{j=1}^{k} [\boldsymbol{\Pi}^{(m)}]_{Y_n,j} )}(V_{n+1}),\\
Y_{0}=& X_0 \; \; \mathrm{a.s.},
\end{split}
\end{equation}
forms the Markov chain $\{Z_n^{(m)}\}$ with transition probability matrix $\boldsymbol{\Pi}^{(m)}$ (\ref{eq:MatrM}).
\end{lem}

\begin{proof}
By applying the Bayes' formula and the law of total probability,
\begin{align}
\mathbb{P}\{Y_{n+1}=i_{n+1} \vert Y_0=i_0\} =&
\sum_{i_1\in\mathbb{S}}\dots  \sum_{i_{n}\in\mathbb{S}} \mathbb{P}\{Y_{n+1}=i_{n+1}, Y_{n}=i_{n},\dots, Y_1=i_1 \vert Y_0=i_0\}  \nonumber\\
=& \sum_{i_1\in\mathbb{S}}\dots  \sum_{i_{n}\in\mathbb{S}}  \mathbb{P}\{Y_1=i_1 \vert Y_0=i_0\} \mathbb{P}\{Y_2=i_2 \vert Y_1=i_1, Y_0=i_0\}   \nonumber\\
&  \hspace{1cm}\times \dots \times \mathbb{P}\{Y_{n}=i_{n}\vert Y_{n-1}=i_{n-1},\dots, Y_0=i_0\}  \nonumber \\
&  \hspace{1cm}\times \mathbb{P}\{Y_{n+1}=i_{n+1}\vert Y_{n}=i_{n},\dots, Y_0=i_0\}. \label{eq:pers1}
\end{align}
On account that $\{V_n\}$ is a series of independent random variables independent of $U_0$, $\mathbb{P}\{Y_{k}=i_{k}\vert Y_{k-1}=i_{k-1},\dots, Y_0=i_0\}=\mathbb{P}\{Y_{k}=i_{k}\vert Y_{k-1}=i_{k-1}\}$. Thus, 
\begin{align}\label{eq:ID2}
 \mathbb{P}\{Y_{k}=i_{k}\vert Y_{k-1}=i_{k-1},\dots,Y_0=i_0\}=\mathbf{e}_{i_{k-1}}^{\top} \boldsymbol{\Pi}^{(m)}  \mathbf{e}_{i_{k}},
\end{align}
for $k=1,\dots,n+1$. As $\sum_{i_k\in\mathbb{S}}\mathbf{e}_{i_k}\mathbf{e}_{i_k}^{\top}=\mathbf{I}_{p\times p}$, we have following (\ref{eq:ID2}) and (\ref{eq:pers1}), 
\begin{align*}
\mathbb{P}\{Y_{n+1}=i_{n+1} \vert Y_0=i_0\} =&
\sum_{i_1\in\mathbb{S}}\dots  \sum_{i_{n}\in\mathbb{S}} \mathbf{e}_{i_0}^{\top} \boldsymbol{\Pi}^{(m)} \mathbf{e}_{i_1}\mathbf{e}_{i_1}^{\top}\boldsymbol{\Pi}^{(m)} \dots  \boldsymbol{\Pi}^{(m)}\mathbf{e}_{i_{n}} \mathbf{e}_{i_{n}}^{\top} [\boldsymbol{\Pi}^{(m)}]  \mathbf{e}_{i_{n+1}} \nonumber \\
=&  \mathbf{e}_{i_0}^{\top} \boldsymbol{\Pi}^{(m)} \Big(\sum_{i_1\in\mathbb{S}} \mathbf{e}_{i_1}\mathbf{e}_{i_1}^{\top} \Big)\boldsymbol{\Pi}^{(m)} \dots \boldsymbol{\Pi}^{(m)} \Big(\sum_{i_n\in\mathbb{S}} \mathbf{e}_{i_{n}} \mathbf{e}_{i_{n}}^{\top} \Big) [\boldsymbol{\Pi}^{(m)}]  \mathbf{e}_{i_{n+1}} \nonumber\\
=& \mathbf{e}_{i_0}^{\top} \big[\boldsymbol{\Pi}^{(m)} \big]^{n+1} \mathbf{e}_{i_{n+1}},   
\end{align*}
showing that $\{Y_n\}$ forms the Markov chains with transition matrix $\boldsymbol{\Pi}^{(m)}$. \exit
\end{proof}

Next, define a Bernoulli random variable $\Phi^{(m)}=\mathbb{1}_{\{\phi=m\}}$ and the function: 
\begin{align}\label{eq:funcF}
F(i,\Phi,V)=\sum_{k=1}^p\sum_{m=1}^M k \Phi^{(m)} \mathbb{1}_{ [\sum_{j=1}^{k-1} [\boldsymbol{\Pi}^{(m)}]_{i,j}, \; \sum_{j=1}^{k} [\boldsymbol{\Pi}^{(m)}]_{i,j} )}(V),
\end{align}
where we set $\sum_{j=1}^{0} [\boldsymbol{\Pi}^{(m)}]_{i,j}=0$ for all $i=1,\dots,p$, $m=1,\dots,M$ and $n\in \mathbb{N}$. Following (\ref{eq:funcF}), consider a finite mixture of Markov chains $\{Z_n\}$ defined by
\begin{equation}\label{eq:chain}
\begin{split}
Z_{n+1}=F(Z_n,\boldsymbol{\Phi},V_{n+1}) \quad \mathrm{with} \quad
Z_0=X_0 \;\; \mathrm{a.s.},
\end{split}
\end{equation}
where $\boldsymbol{\Phi}$ denotes $(1\times M)-$vector of speed regime, i.e., $\boldsymbol{\Phi}=(\Phi^{(1)},\dots,\Phi^{(M)}).$

\begin{prop} The representation (\ref{eq:chain}) yields the mixture of Markov chains
\begin{align}\label{eq:chainmixture1}
Z_n=\sum_{m=1}^M \Phi^{(m)} Z_n^{(m)}, \quad \textrm{with\; $Z_0=X_0$},
\end{align}
where $Z_{n}^{(m)}$ satisfies the recursive equation $Z_{n+1}=F(Z_n,\Phi^{(m)}=1,V_{n+1})$, see eqn.  (\ref{eq:chainmixture2}), while the $n$-step transition probability matrix of $\{Z_n\}$ is defined by
\begin{align}\label{eq:stochmatrix}
\mathbf{P}^{(n)} = \sum_{m=1}^M \mathbf{S}^{(m)} \big[\boldsymbol{\Pi}^{(m)}\big]^n.
\end{align}
\end{prop}
\begin{proof}
The representation (\ref{eq:chainmixture1}) follows from (\ref{eq:funcF}) and (\ref{eq:chain})  and applying the Fubini's principle to (\ref{eq:funcF}). By the Bayes formula and law of total probability, 
\begin{align*}
\big[\mathbf{P}^{(n)}\big]_{i,j}=&\mathbb{P}\{Z_n=j\vert Z_0=i\}\\
=&\sum_{m=1}^M \mathbb{P}\{\Phi^{(m)}=1\vert Z_0=i\}\mathbb{P}\{Z_n=j\vert \Phi^{(m)}=1, Z_0=i\}\\
=&\sum_{m=1}^M  \mathbb{P}\{\Phi^{(m)}=1\vert Z_0=i\} \mathbb{P}\{Z_n^{(m)}=j\vert Z_0^{(m)}=i\} \\
=& \sum_{m=1}^M s_i^{(m)} \big[\big[\boldsymbol{\Pi}^{(m)}\big]^n]_{i,j}=\mathbf{e}_i^{\top} \sum_{m=1}^M \mathbf{S}^{(m)} \big[\boldsymbol{\Pi}^{(m)}\big]^n \mathbf{e}_j,
\end{align*}
where the last equality follows from taking account of Lemma \ref{lem:mainlemma}, leading to the establishment of identity (\ref{eq:stochmatrix}), given $\mathbf{S}^{(m)}$ is a $(p\times p)-$diagonal matrix. \exit
\end{proof}
\begin{Rem}
It is worth noticing that if the underlying Markov chains $\{Z_n^{(m)}:m=1,\dots,M\}$ have the same transition probability matrix $\boldsymbol{\Pi}$, which is the case for the mixture model \cite{Frydman2005}, one can show that $\{Z_{n}\}$ (\ref{eq:chain}) has the same distribution as $\{Y_{n}\}$ (\ref{eq:chainmixture2}). The latter is used in \cite{Resnick} to generate the Markov chains $\{Y_n\}$. 
\end{Rem}

\begin{figure}[htp!]
 \includegraphics[scale=.85,center]{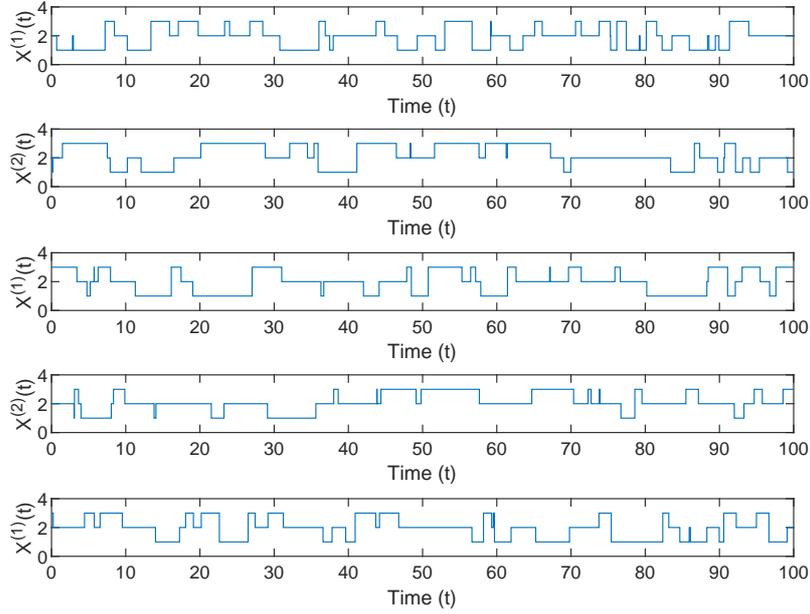}
 \caption{Sample paths of the mixture of Markov jump process $X$ (\ref{eq:mixture})}
 \label{fig:Fig1}
\end{figure}

\subsection{Finite mixture of Markov jump processes}

The epoch time $\{T_n\}$ of the mixture process $X$ (\ref{eq:mixture}) is defined by
\begin{align}\label{eq:epoch}
T_{n+1}=T_n- \sum_{k=1}^M \Phi^{(k)} \frac{\log W_n}{q_{Z_n}^{(k)}}, \;\; \textrm{with $T_0=0$ a.s.}
\end{align}
\begin{theo}\label{theo:maintheo}
Define a continuous-time stochastic process $X=\{X(t):t\geq 0\}$,
\begin{align}\label{eq:simulated}
X(t)=\sum_{n=0}^{\infty} Z_n \mathbb{1}_{[ T_n, T_{n+1})}(t) \;\; \mathrm{with} \;\; X(0)=X_0.
\end{align}
Then, the stochastic process $X$ has the following properties:
\begin{enumerate}
\item[(i)] it has the representation of the mixture process (\ref{eq:mixture}) with
\begin{align}\label{eq:processXmt}
X^{(m)}(t):=\sum_{n=0}^{\infty} Z_n^{(m)}\mathbb{1}_{[T_n, T_{n+1})}(t),
\end{align}
\item[(ii)]  and has the same distribution (\ref{eq:transM}) as the mixture process $X$ (\ref{eq:mixture}).
\end{enumerate}
\end{theo}

Figure \ref{fig:Fig1} displays Monte Carlo simulation of the sample paths of $X$ (\ref{eq:simulated}). As we can see, at every fixed point of time $t\geq 0$ and state $i\in\mathbb{S}$ the observed paths is comprised of a mixture of two Markov jump processes $X^{(1)}(t)$ and $X^{(2)}(t)$.

In order to establish the theorem, the following results are required.

\begin{lem}\label{lem:reqlem1}
Define the following transition probability matrix $\mathbf{P}_{i,j}^{(m)}(t)$:
\begin{align}\label{eq:transMatrix}
\mathbf{P}_{i,j}^{(m)}(t):=\mathbb{P}\{X(t)=j\vert \Phi^{(m)}=1,X(0)=i\}.
\end{align}
Then, for all $t\geq 0$, the function $t\rightarrow \mathbf{P}_{i,j}^{(m)}(t)$ solves the integral equation: 
\begin{align}\label{eq:inteq}
\mathbf{P}_{i,j}^{(m)}(t)=e^{-q_i^{(m)} t} \delta_{ij}  + \int_0^t  q_i^{(m)} e^{-q_i^{(m)} u} \sum_{l\neq i} \pi_{i,l}^{(m)}  \mathbf{P}_{l,j}^{(m)}(t-u) du.
\end{align}
\end{lem}

\begin{proof}
Since $T_1$ is the first jump time of $X$, the number of steps $n$ until the Markov chain $Z^{(m)}$ makes a jump from state $i$ to state $j\neq i$ has geometric distribution with probability $\mathbf{e}_i^{\top}[1-\boldsymbol{\Pi}^{(m)}]^{n-1} \boldsymbol{\Pi}^{(m)}\mathbf{e}_j$. By independence of $W_0$, 
\begin{align*}
 &\mathbb{P}\{X(t)=j, T_1>t \vert \Phi^{(m)}=1, X(0)=i\}\\
 &\hspace{2cm}= 
 \sum_{n=1}^{\infty} \mathbb{P}\{Z_n=j, T_1>t\vert \Phi^{(m)}=1, Z_0=i\}\\
 &\hspace{2cm}=\sum_{n=1}^{\infty} \mathbb{P}\Big\{Z_n^{(m)}=j, -\frac{\log W_0}{q_i^{(m)}}>t \big\vert Z_0^{(m)}=i \Big\}\\
  &\hspace{2cm}=\sum_{n=1}^{\infty}  e^{-q_i^{(m)}t} \mathbf{e}_i^{\top}(1-\boldsymbol{\Pi}^{(m)})^{n-1} \boldsymbol{\Pi}^{(m)}\mathbf{e}_j \\
%  &\hspace{2cm}= e^{-q_i^{(m)}t} \mathbf{e}_i^{\top}\Big(\sum_{n=1}^{\infty}[1-\boldsymbol{\Pi}^{(m)}]^{n-1}\Big) \boldsymbol{\Pi}^{(m)} \mathbf{e}_j\\
  &\hspace{2cm}=e^{-q_i^{(m)}t} \delta_{ij},
 \end{align*}
 with $\delta_{ij}= \mathbf{e}_i^{\top} \mathbf{e}_j$, provided the transition probability matrix $\boldsymbol{\Pi}^{(m)}$ is invertible. 
 
Furthermore, following (\ref{eq:epoch}) the epoch time $T_n$ of $X$ can be singled out as  
\begin{align}\label{eq:epoch2}
T_n=-\sum_{k=0}^{n-1}\sum_{m=1}^M \Phi^{(m)} \frac{\log W_k}{q_{Z_k}^{(m)}} \quad \textrm{with \; $T_0=0$ \; a.s.}
\end{align}
Again, given that $T_1$ is the first jump time of $X$ moving from state $i$ to $j\neq i$, we have by independence of $W_0$, Bayes' formula and Fubini's principle that
\begin{align*}
&\mathbb{P}\{X(t)=j, T_1\leq t \vert \Phi^{(m)}=1, X(0)=i\}\\
&\hspace{0cm}=\sum_{l\neq i} \mathbb{P}\Big\{X(t)=j, -\frac{\log W_0}{q_i^{(m)}} \leq t, Z_1= l \big\vert \Phi^{(m)}=1, Z_0=i \Big\}\\
&\hspace{0cm}=\sum_{l\neq i} \int_0^t \mathbb{P}\Big\{  -\frac{\log W_0}{q_i^{(m)}}\in du \big \vert \Phi^{(m)}=1, Z_0=i\Big\}  \mathbb{P}\big\{Z_1=l\vert \Phi^{(m)}=1, Z_0=i \big\}\\
&\hspace{3cm} \times \mathbb{P}\Big\{ X(t)=j \big\vert \Phi^{(m)}=1,   -\frac{\log W_0}{q_i^{(m)}}=u, Z_1=l, Z_0=i \Big\}\\
&\hspace{0cm}=\int_0^t q_i^{(m)}  e^{-q_i^{(m)}u} \sum_{l\neq i} \pi_{i,l}^{(m)} \mathbb{P}\Big\{ X(t)=j \big\vert \Phi^{(m)}=1,   -\frac{\log W_0}{q_i^{(m)}}=u, Z_1=l, Z_0=i \Big\}du.
\end{align*}
The proof is accomplished once we have shown that the following identity holds:
\begin{align}\label{eq:persID}
\mathbb{P}\Big\{ X(t)=j \big\vert \Phi^{(m)}=1,   -\frac{\log W_0}{q_i^{(m)}}=u, Z_1=l, Z_0=i \Big\}=\mathbf{P}_{l,j}^{(m)}(t-u).
\end{align}
To this end, recall that the conditional probability in (\ref{eq:persID}) can be simplified as
\begin{align*}
&\mathbb{P}\Big\{ X(t)=j \big\vert \Phi^{(m)}=1,   -\frac{\log W_0}{q_i^{(m)}}=u, Z_1=l, Z_0=i\Big\}\\
&=\sum_{n=1}^{\infty} \mathbb{P}\big\{Z_n=j, T_n\leq t < T_{n+1}\big \vert \Phi^{(m)}=1,  -\frac{\log W_0}{q_i^{(m)}}=u, Z_1=l,Z_0=i\big\}\\
&=\sum_{n=1}^{\infty} \mathbb{P}\Big\{Z_n=j, -\sum_{k=1}^{n-1} \frac{\log (W_k)}{q_{Z_k}^{(m)}}\leq t-u < -\sum_{k=1}^{n} \frac{\log (W_k)}{q_{Z_k}^{(m)}} \Big\vert \Phi^{(m)}=1, Z_1=l, Z_0=i \Big\}\\
&= \sum_{n=1}^{\infty} \mathbb{P}\Big\{Z_{n-1}=j, -\sum_{k=1}^{n-1} \frac{\log (W_k)}{q_{Z_{k-1}}^{(m)}}\leq t-u< -\sum_{k=1}^{n} \frac{\log (W_k)}{q_{Z_{k-1}}^{(m)}} \Big\vert \Phi^{(m)}=1, Z_0=l \Big\},
\end{align*}
where in the last equality we have used the fact that conditional on $\Phi^{(m)}=1$, the discrete-time mixture process $Z_n$ is a Markov chain moving according to $Z_n^{(m)}$, which by Lemma \ref{lem:mainlemma} has the memoryless property. Next, define new indexes $n^{\prime}=n-1$, $k^{\prime}=k-1$, and $W_k^{\prime}=W_{k^{\prime}+1}$. By doing so, we can rewrite
\begin{align*}
& \sum_{n=1}^{\infty} \mathbb{P}\Big\{Z_{n-1}=j, -\sum_{k=1}^{n-1} \frac{\log (W_k)}{q_{Z_{k-1}}^{(m)}}\leq t-u< -\sum_{k=1}^{n} \frac{\log (W_k)}{q_{Z_{k-1}}^{(m)}} \Big\vert \Phi^{(m)}=1, Z_0=l \Big\}\\
 &= \sum_{n^{\prime}=0}^{\infty} \mathbb{P}\Big\{Z_{n^{\prime}}=j, -\sum_{k^{\prime}=0}^{n^{\prime}-1} \frac{\log (W_k^{\prime})}{q_{Z_{k^{\prime}}}^{(m)}}\leq t-u< -\sum_{k^{\prime}=0}^{n^{\prime}} \frac{\log (W_k^{\prime})}{q_{Z_{k^{\prime}}}^{(m)}} \Big\vert \Phi^{(m)}=1, Z_0=l \Big\}\\
 &=\mathbf{P}_{l,j}^{(m)}(t-u).
 \end{align*}

The integral equation (\ref{eq:inteq}) is obtained by invoking the law of total probability and the Bayes' formula for conditional probability to (\ref{eq:transMatrix}) to get
\begin{align*}
\mathbf{P}_{i,j}^{(m)}(t)=& \mathbb{P}\{X(t)=j, T_1>t \vert \Phi^{(m)}=1, X(0)=i\} \\
&+ \mathbb{P}\{X(t)=j, T_1\leq t \vert \Phi^{(m)}=1, X(0)=i\}. 
\end{align*}
The final result (\ref{eq:inteq}) is obtained by collecting the two pieces together. \exit
\end{proof}

\medskip

\begin{prop}\label{prop:mainprop}
The integral equation (\ref{eq:inteq}) has an explicit solution:
\begin{align}\label{eq:probmat}
\mathbf{P}_{i,j}^{(m)}(t)=\mathbf{e}_i^{\top} e^{\mathbf{Q}^{(m)}t} \mathbf{e}_j, \quad \textrm{i.e.,} \quad  \mathbf{P}^{(m)}(t)=e^{\mathbf{Q}^{(m)}t}.
\end{align}
%
%where $\mathbf{Q}^{(m)}$ is a $(p\times p)-$intensity matrix defined by 
%%
%\begin{align*}
%[\mathbf{Q}^{(m)}]_{i,j}=
%\begin{cases}
%-q_i^{(m)}, &\mathrm{for} \; j=i \\
%q_i^{(m)} [\boldsymbol{\Pi}^{(m)}]_{i,j}, &\mathrm{for} \; j\neq i \\
%\end{cases}
%\end{align*}
%%
\end{prop}
\begin{proof}
On account of the fact that $\sup_t \mathbf{P}_{i,j}^{(m)}(t)\leq 1$ and $\sum_{l\neq i} \pi_{i,l}^{(m)}=1$, the function inside the integral in (\ref{eq:inteq}) is uniformly bounded. Thus, the integral is a continuous function of $t$, and therefore so is the function $\mathbf{P}_{i,j}^{(m)}(t)$. Hence, as a result, $\mathbf{P}_{i,j}^{(m)}(t)$ is absolutely continuous w.r.t. to Lebesgue measure $dt$, i.e., it is continuously differentiable. Applying change of variable $s=t-u$, we obtain
\begin{align*}
\mathbf{P}_{i,j}^{(m)}(t)=e^{-q_i^{(m)}t} \Big(\delta_{ij} + q_i^{(m)}\int_0^t  e^{q_i^{(m)}s} \sum_{l\neq i} \pi_{i,l}^{(m)} \mathbf{P}_{l,j}^{(m)}(s) ds\Big).
\end{align*}
As the function $\mathbf{P}_{i,j}^{(m)}(t)$ is continuously differentiable, we have 
\begin{align*}
\frac{d}{dt}\mathbf{P}_{i,j}^{(m)}(t)=& - q_i^{(m)} e^{-q_i^{(m)}t}\Big(\delta_{ij} + q_i^{(m)}\int_0^t  e^{q_i^{(m)}s} \sum_{l\neq i} \pi_{i,l}^{(m)} \mathbf{P}_{l,j}^{(m)}(s) ds\Big)\\
& +e^{-q_i^{(m)}t} \Big(q_i^{(m)} e^{q_i^{(m)}t}  \sum_{l\neq i} \pi_{i,l}^{(m)} \mathbf{P}_{l,j}^{(m)}(t)\Big) \\
=& - q_i^{(m)} \mathbf{P}_{i,j}^{(m)}(t) + q_i^{(m)} \sum_{l\neq i} \pi_{i,l}^{(m)}\mathbf{P}_{l,j}^{(m)}(t) \\
=&\sum_{l\in\mathbb{S}} \mathbf{Q}_{i,l}^{(m)}\mathbf{P}_{l,j}^{(m)}(t).
\end{align*}
As $\mathbf{P}^{(m)}(0)=\mathbf{I}$, it is straightforward to check that the linear systems of equation has the unique solution (\ref{eq:probmat}) for the transition matrix $\mathbf{P}^{(m)}(t)$ for $t\geq0$. \exit
\end{proof}

\subsubsection*{Proof of Theorem \ref{theo:maintheo}}
The claim is established by the law of total probability and Bayes' formula for conditional probability taking account of Lemma  \ref{lem:reqlem1} and Proposition \ref{prop:mainprop}, i.e.,
\begin{align*}
[\mathbf{P}(t)]_{i,j}=&\mathbb{P}\{X(t)=j \vert X(0)=i\}\\=&\sum_{m=1}^M \mathbb{P}\{\Phi^{(m)}=1\vert X(0)=i\}\ \mathbb{P}\{X(t)=j \vert \Phi^{(m)}=1, X(0)=i\}\\
=& \sum_{m=1}^M s_i^{(m)} \mathbf{P}_{i,j}^{(m)}(t) 
= \mathbf{e}_i^{\top}\sum_{m=1}^M \mathbf{S}^{(m)} \mathbf{P}^{(m)}(t) \mathbf{e}_j.
%=& \mathbf{e}_i^{\top}\sum_{m=1}^M \mathbf{S}^{(m)} e^{\mathbf{Q}^{(m)}(t)} \mathbf{e}_j.
\end{align*}
Inserting the expression of $Z_n$ (\ref{eq:chainmixture1}) in (\ref{eq:simulated}) yields the mixture (\ref{eq:mixture}). \exit

%\medskip
%
%\textbf{Pseudo algorithm for the simulation}
%
%\begin{enumerate}
%
%\item[(1)] \textbf{Step 1} Choose initial state $X_0=i_0$ with the probability $\boldsymbol{\pi}$ using (\ref{eq:initialcond}). 
%
%\item[(2)] \textbf{Step 2} On the chosen initial state $i_0$, select the speed regime $\phi=m$ with the probability $\mathbf{s}_{i_0}=(s_{i_0}^{(1)},\dots,s_{i_0}^{(M)})$ using (\ref{eq:initialregime}). Then, set $\Phi^{(m)}=1$.
%
%\item[(3)] \textbf{Step 3} Make transition from time zero to a specified final time $T>0$ according to the selected Markov jump process $\{X_t^{(m)}\}_{t\geq 0}$ starting from the initial state $i_0$ by applying the recursive equations (\ref{eq:chain})-(\ref{eq:simulated}). 
%
%\end{enumerate}
%

\section{Estimation with complete information}
Statistical estimation of $\{q_{ij}^{(m)}\}$, $\{q_{i}^{(m)}\}$, and $\{s_i^{(m)}\}$ was discussed in \cite{Frydman2005} for special structure of intensity matrix $\mathbf{Q}^{(m)}$, where it is assumed that $\mathbf{Q}^{(m)}=\boldsymbol{\Psi}^{(m)}\mathbf{Q}$ with $\boldsymbol{\Psi}^{(m)}=\mathrm{diag}(\psi_1^{(m)},\dots,\psi_p^{(m)})$ for $1\leq m\leq M-1$ and $\boldsymbol{\Psi}^{(M)}=\mathbf{I},$ i.e., $q_{ij}^{(m)}=\psi_i^{(m)} q_{ij}$ implying that $\pi_{ij}^{(m)}=\pi_{ij}$ for all $m=1,\dots,M$, see (\ref{eq:MatrM}). This is equivalent to imposing the condition on the embedded Markov chains $Z_n^{(m)}$ to have the same transition probability matrix $\boldsymbol{\Pi}$, i.e., $\boldsymbol{\Pi}^{(m)}=\boldsymbol{\Pi}$ for all $m$.

This paper attempts to generalize the estimation method \cite{Frydman2005} for inferring the distribution parameters $\{q_{ij}^{(m)}\}$, $\{q_{i}^{(m)}\}$, and $\{s_i^{(m)}\}$ of $X$ for a general structure of $\mathbf{Q}^{(m)}$ as well as to estimate $\{\pi_i\}$, which was not discussed in \cite{Frydman2005}. Importantly, as opposed to the EM estimator $\widehat{s}_i^{(m)}$ given in \cite{Frydman2005}, our estimate for $s_i^{(m)}$ sums to one, i.e., $\sum_{m=1}^M \widehat{s}_i^{(m)} =1$, for $i=1,\dots,p$, which is due to the constraint (\ref{eq:portion}). The results can be used to recover the estimation for the restricted mixture \cite{Frydman2005}.

\subsection{Maximum likelihood estimation}

This section discusses maximum likelihood estimation of the distribution of $X$ based on complete observations of $X$, where the underlying process driving the evolution of $X$ is known. To start with, we assume that $N$ independent realizations $\{X_k\}$ of $X$ are observed continuously on the time interval $[0,T]$ with $0<T<\infty.$ For notational convenience, we use the following conventions:
\begin{equation}\label{eq:data}
\begin{split}
\Phi_k^{(m)}=&\mathbb{1}_{\{X_k=X^{(m)}\}}\\[6pt]
B_i^{(k)}=&\mathbb{1}_{\{X_k(0)=i\}}\\[6pt]
N_{ij}^{(k)}=&\sum_{l=1}^{N} \mathbf{1}_{\{ X_k(lh)=j, X_k((l-1)h)=i\}}, \; \mathrm{with} \; \;h=T/N, \; N\in \mathbb{N} \\[6pt]
N_i^{(k)}=& \sum_{j\neq i} N_{ij}^{(k)}\\[6pt]
Z_i^{(k)}=&\int_0^{T}\mathbb{1}_{\{X_k(u)=i\}}du,
\end{split}
\end{equation}
for $i,j=1,\dots,p$, $k=1,\dots,N$ and $m=1,\dots,M$. More precisely, $N_{ij}^{(k)}$ counts the (number) of transitions of $X_k$ from state $i$ to state $j\neq i$, $N_{i}^{(k)}$ counts the number of transitions from state $i$, whereas $Z_i^{(k)}$ represents the occupation time of $X_k$  in state $i$ when $X_k$ is observed on the interval $[0,T]$, while $B_i^{(k)}$ counts the number of realizations of $X_k$ starting in state $i$ at time zero. Also, note that
\begin{align}\label{eq:cond0}
\sum_{m=1}^M \Phi_k^{(m)}=1  \quad \mathrm{and} \quad \sum_{i=1}^p B_i^{(k)}=1, \quad \textrm{for $k=1,\dots,N$},
\end{align}
which in turn implies that  
\begin{align}\label{eq:cond1}
\sum_{i=1}^p \sum_{k=1}^N \sum_{m=1}^M \Phi_k^{(m)} B_i^{(k)}=N.
\end{align}

Under complete information, the statistics (\ref{eq:data}) are assumed to be available for the maximum likelihood estimation of the distribution parameters of $X$ (\ref{eq:mixture}). In the sequel below we denote by $\boldsymbol{\theta}=(\boldsymbol{\pi},\mathbf{Q}^{(m)},\mathbf{S}^{(m)})$ and $f_{\boldsymbol{\theta}}(X_k,\Phi_k^{(m)})$ the joint probability density function of the observations $X_k$ and $\Phi_k^{(m)}$, indicator variable which provides information on which underlying process that drives $X_k$. 

Suppose that $X_k$ chooses its initial state $i_k$ to start with randomly at probability $\pi_{i_k}$. On account that the bivariate process $(X_k,\Phi_k^{(m)})$ is Markovian, it follows from applying the Bayes' formula for conditional probability that
\begin{align}
f_{\boldsymbol{\theta}}(X_k,\Phi_k^{(m)})=&f_{\boldsymbol{\theta}}(X_k(0)=i_k)f_{\boldsymbol{\theta}}\big(\Phi_k^{(m)} \vert X_k(0)=i_k\big) \nonumber \\ &\times f_{\boldsymbol{\theta}}\big(X_k \vert \Phi_k^{(m)}, X_k(0)=i_k\big)  \nonumber\\[6pt]
=& \big( \pi_{i_k} s_{i_k}^{(m)} \big)^{B_{i_k}^{(k)}} f_{\boldsymbol{\theta}}\big(X_k \vert \Phi_k^{(m)}, X_k(0)=i_k\big). \label{eq:pdf1}
\end{align}
Conditional on knowing $\Phi_k^{(m)}$, $f_{\boldsymbol{\theta}}\big(X_k \vert \Phi_k^{(m)}, X_k(0)=i_k\big)$ represents the likelihood function of observing the sample paths of $X_k$ under the Markov process $X^{(m)}$ for which the likelihood is given by (e.g. Albert \cite{Albert}, Basawa and Rao \cite{Basawa}):
\begin{align*}
f_{\boldsymbol{\theta}}\big(X_k \vert \Phi_k^{(m)}, X_k(0)=i_k\big)
 =\prod_{i=1}^p \prod_{j\neq i}^p \big(q_{ij}^{(m)}\big)^{N_{ij}^{(k)}} e^{-q_{ij}^{(m)} Z_i^{(k)}}.
\end{align*}
Hence, following (\ref{eq:pdf1}), the likelihood contribution $L_k^{(m)}$ of $(X_k,\Phi_k^{(m)})$ is given by
\begin{align*}
L_k^{(m)}=& \prod_{i=1}^p (s_i^{(m)} \pi_i)^{B_i^{(k)}}\prod_{i=1}^p \prod_{j\neq i}^p \big(q_{ij}^{(m)}\big)^{N_{ij}^{(k)}} e^{-q_{ij}^{(m)} Z_i^{(k)}}.
\end{align*}

The likelihood contribution of all realizations $\{X_k\}$ is therefore given by  
\begin{align}\label{eq:likelihood}
L= \prod_{k=1}^N \prod_{m=1}^M  \big[ L_k^{(m)} \big]^{\Phi_k^{(m)}}.
\end{align}

Notice that the likelihood (\ref{eq:likelihood}) reduces to the one given in \cite{Albert} and  \cite{Basawa} when $s_i^{(m)}=1$ for all $i\in\mathbb{S}$ and $m$, and $M=1$, which implies $\Phi_k^{(m)}=1$ and $q_{ij}^{(m)}=q_{ij}$ for all $i,j\in\mathbb{S}$. Equivalently in terms of the log-likelihood function, we have 
\begin{align}\label{eq:log-like}
\log L=&\sum_{k=1}^N \sum_{m=1}^M \Phi_k^{(m)} \log L_k^{(m)} \nonumber\\
=& \sum_{k=1}^N \sum_{m=1}^M \Phi_k^{(m)} \Big[ \sum_{i=1}^p B_i^{(k)} \log \big( s_i^{(m)} \pi_i\big)  \\
& \hspace{2.5cm} + \sum_{i=1}^p \sum_{j\neq i}^p N_{ij}^{(k)} \log q_{ij}^{(m)} - \sum_{i=1}^p \sum_{j\neq i}^p q_{ij}^{(m)} Z_i^{(k)} \Big]. \nonumber
\end{align}
Recall that the probability $\{\pi_i\}$ and $\{s_i^{(m)}\}$ satisfy the constraint, see (\ref{eq:portion}):
\begin{align}\label{eq:cond2}
\sum_i^p \pi_i=1 \quad \mathrm{and} \quad \sum_{m=1}^M s_i^{(m)}=1, \;\; \textrm{for all $i=1,\dots,p$}.
\end{align}

The maximum likelihood estimators for $\pi_i$, $q_{ij}^{(m)}$, $q_i^{(m)}$, and $s_i^{(m)}$ are given explicitly in terms of the statistics (\ref{eq:data}). The results are summarized below.

\begin{theo}
The maximum likelihood estimates of $\pi_i$, $q_{ij}^{(m)}$, $q_i^{(m)}$, and $s_i^{(m)}$ are 
\begin{eqnarray}
\widehat{\pi}_{N,i} &=&\frac{1}{N}\sum_{k=1}^N B_i^{(k)}, \label{eq:est1}\\[6pt]
\widehat{q}_{N,ij}^{(m)}&=&\;\frac{\sum_{k=1}^N \Phi_k^{(m)} N_{ij}^{(k)}}{ \sum_{k=1}^N \Phi_k^{(m)} Z_{i}^{(k)}},\label{eq:est2}\\[6pt]
\widehat{q}_{N,i}^{(m)} &=&\frac{\sum_{k=1}^N \Phi_k^{(m)} N_i^{(k)}}{\sum_{k=1}^N \Phi_k^{(m)} Z_{i}^{(k)}}, \label{eq:est3} \\[6pt]
\widehat{s}_{N,i}^{(m)} &=&\frac{\sum_{k=1}^N \Phi_k^{(m)} B_i^{(k)}}{\sum_{k=1}^N  B_i^{(k)}}. \label{eq:est4}
\end{eqnarray}
 \end{theo}
\begin{proof}
To find the estimators $\widehat{q}_{N,ij}^{(m)}$, $\widehat{q}_{N,i}^{(m)}$, $\widehat{s}_{N,i}^{(m)}$ and $\widehat{\pi}_{N,i}$ of the parameters distribution of the mixture process $X$, we introduce the Lagrangian function:
\begin{align*}
\mathcal{L}= \log L- \lambda \Big(\sum_{i=1}^p \pi_i -1\Big) - \sum_{i=1}^p \gamma_i \Big(\sum_{m=1}^M s_i^{(m)}-1\Big),
\end{align*}
where $\lambda$ and $\gamma_i$ are the corresponding Lagrange multipliers of the constraint (\ref{eq:cond2}). Applying the first order Euler condition w.r.t $\pi_i$ to the Lagrangian $\mathcal{L}$, 
\begin{align*}
\frac{\partial \mathcal{L}}{\partial \pi_i}=\sum_{k=1}^N \sum_{m=1}^M \Phi_k^{(m)}\frac{B_i^{(k)}}{\pi_i}-\lambda=0 \; \Longrightarrow \lambda=N,
\end{align*}
taking the note on  (\ref{eq:cond0}), (\ref{eq:cond1}) and (\ref{eq:cond2}). The estimator $\widehat{\pi}_{N,i}$ is given by (\ref{eq:est1}). 

To get the estimator $\widehat{q}_{N,ij}^{(m)}$, we set the following Euler equation:
\begin{align*}
\frac{\partial \mathcal{L}}{\partial q_{ij}^{(m)}}=\sum_{k=1}^N \Phi_k^{(m)} \Big(\frac{N_{ij}^{(k)}}{q_{ij}^{(m)}} - Z_i^{(k)}\Big)=0,
\end{align*}
solving which for $q_{ij}^{(m)}$ gives the estimator $\widehat{q}_{N,ij}^{(m)}$ defined by (\ref{eq:est2}).

\medskip

Given that $q_i^{(m)}$ satisfies (\ref{eq:matq}), we have following (\ref{eq:data}), $\widehat{q}_{N,i}^{(m)}=\sum\limits_{j\neq i} \widehat{q}_{N,ij}^{(m)}.$

Finally, we set the following Euler equation for $s_i^{(m)}$:
\begin{align*}
\frac{\partial \mathcal{L}}{\partial s_i^{(m)}}=\sum_{k=1}^N \Phi_k^{(m)} \frac{B_i^{(k)}}{s_i^{(m)}} -\gamma_i=0 \; \Longrightarrow \gamma_i=\sum_{k=1}^N B_i^{(k)},
\end{align*}
on account of (\ref{eq:cond0}) and (\ref{eq:cond2}). Hence, the estimator $\widehat{s}_i^{(m)}$ is given by (\ref{eq:est4}). \exit
\end{proof}

\subsubsection{Restricted mixture of Markov jump processes}\label{sec:restrictedmixture}

The mixture model and its EM estimation were first discussed in \cite{Frydman2005} for a special class of mixture process in which case the element $q_{ij}^{(m)}$ of the intensity matrix $\mathbf{Q}^{(m)}$ is defined by $q_{ij}^{(m)}=\psi_i^{(m)} q_{ij}$ implying that each underlying Markov jump process has the same probability of leaving a state, i.e., $\pi_{ij}^{(m)}=\pi_{ij}$. For simplicity, we set following \cite{Frydman2005} $\psi^{(M)}=1$. As a result, the maximum likelihood estimators of $q_{ij}$ and $\psi_i^{(m)}$, for $m=1,\dots,M-1$, $j\neq i$, are given following (\ref{eq:log-like}) and (\ref{eq:cond2}) by 
\begin{align}
\widehat{q}_i=&\frac{\sum_{k=1}^N \Phi_k^{(M)} N_i^{(k)}}{\sum_{k=1}^N \Phi_k^{(M)} Z_i^{(k)}},\label{eq:est5}\\[6pt]
\widehat{\psi}_i^{(m)}=& \frac{\sum_{k=1}^N \Phi_k^{(m)} N_i^{(k)}}{\widehat{q}_i \sum_{k=1}^N \Phi_k^{(m)} Z_i^{(k)}}, \label{eq:est6}\\[6pt]
\widehat{\pi}_{ij}=&\frac{\sum_{k=1}^N N_{ij}^{(k)}}{\sum_{k=1}^N N_i^{(k)}}, \nonumber\\[6pt]
\widehat{q}_{ij}=&\widehat{\pi}_{ij}\widehat{q}_i,  \nonumber\\[6pt]
\widehat{q}_{ij}^{(m)}=&\widehat{\psi}_{i}^{(m)}\widehat{q}_{ij},\label{eq:restestm}
\end{align}
while the estimators $\widehat{s}_i^{(m)}$ and $\widehat{\pi}_i$ are the same as (\ref{eq:est1}) and (\ref{eq:est4}), respectively.

\begin{Rem}\label{rem:restestm}
It is straightforward to see following (\ref{eq:restestm}) that $\widehat{\pi}_{ij}^{(m)}=\widehat{\pi}_{ij}.$
\end{Rem}

\subsection{Consistency of the MLE estimators}

To establish consistency of the MLE estimators (\ref{eq:est1})-(\ref{eq:est4}), the following results are required. For convenience, we write $\Phi^{(m)}:=\Phi_1^{(m)}$,$N_{ij}:=N_{ij}^{(1)}$ and $Z_i:=Z_i^{(1)}$.

\begin{lem}\label{lem:lem1}
For given $i,j=1,\dots,p$ and $m=1,\dots,M$, we have
\begin{eqnarray}
\mathbb{E}\big\{\Phi^{(m)}N_{ij}\big\}&=&q_{ij}^{(m)} \int_0^T \mathbb{P}\{X(u)=i, \Phi^{(m)}=1\}du \label{eq:Id1},\\[5pt] 
\mathbb{E}\big\{\Phi^{(m)}Z_i\big\}&=&  \int_0^T \mathbb{P}\{X(u)=i, \Phi^{(m)}=1\}du. \label{eq:Id2}
\end{eqnarray}
\end{lem}
\begin{proof}
Recall that $\Phi^{(m)}N_{ij}$ counts the number of transition of $X$ (\ref{eq:mixture}), which is driven by the underlying Markov process $X^{(m)}$, between state $i$ and $j$ over the period of time $[0,T]$. The proof of (\ref{eq:Id1}) can follow similar approach to the proof of Theorem 5.1(a) in \cite{Albert}. To be more precise, since $\Phi^{(m)}N_{ij}=\sum_{k=1}^{N} \mathbf{1}_{\{ X(kh)=j, X((k-1)h)=i, \Phi^{(m)}=1\}}$, then by the Bayes' formula we obtain,
\begin{align*}
\mathbb{E}\{ \Phi^{(m)}N_{ij} \} =&\sum_{k=1}^{N} \mathbb{P}\{X(kh)=j, X((k-1)h)=i, \Phi^{(m)}=1\}\\
&\hspace{-2cm}=\sum_{k=1}^{N} \mathbb{P}\{X(kh)=j \big \vert X((k-1)h)=i, \Phi^{(m)}=1\} \mathbb{P}\{ X((k-1)h)=i, \Phi^{(m)}=1\}\\
=& q_{ij}^{(m)} \sum_{k=1}^{N}  \mathbb{P}\{ X((k-1)h)=i, \Phi^{(m)}=1\} h + o(h)\\[5pt]
 & \hspace{0.85cm}\xrightarrow[]{ N\rightarrow \infty}  q_{ij}^{(m)}  \int_0^T \mathbb{P}\{ X(u) =i, \Phi^{(m)}=1\} du,
\end{align*}
where the limit is due to Lebesgue dominated convergence theorem. (\ref{eq:Id2}) follows given that $\mathbb{E}\big\{\int_0^T \vert \mathbb{1}_{\{X(u)=i, \Phi^{(m)}=1\}} \vert du\big\}\leq T<\infty$, by which the claim follows from applying Fubini's theorem to the expectation $\mathbb{E}\big\{\int_0^T \mathbb{1}_{\{X(u)=i, \Phi^{(m)}=1\}}du\big\}.$ \exit
\end{proof}

\medskip

Thanks to the results of Lemma \ref{lem:lem1}, consistency of the estimators follows.

\begin{theo}
If for $i=1,\dots,p$ and $t\geq 0$, $\mathbb{P}\{X(t)=i\}>0$, then
\begin{align*}
\lim\limits_{N\rightarrow \infty} \widehat{q}_{N,ij}^{(m)}=&q_{ij}^{(m)}, \quad \quad \lim\limits_{N\rightarrow \infty} \widehat{q}_{N,i}^{(m)}=q_{i}^{(m)}, \\[5pt] \lim\limits_{N\rightarrow \infty} \widehat{s}_{N,i}^{(m)}=&s_{i}^{(m)}, \quad \quad \lim\limits_{N\rightarrow \infty} \widehat{\pi}_{N,i}=\pi_i.
\end{align*}
with probability one, for $i,j=i,\dots,p$ and $m=1,\dots,M$. 

\end{theo}

\begin{proof}
The proof is based on applying the law of large numbers and continuous mapping theorem applied to independent paired observations $(X_k,\Phi_k^{(m)})$. 

For convenience, we write $B_i:=B_i^{(1)}$ and $N_i:=N_i^{(1)}$. Recall that
\begin{align*}
\mathbb{E}\{B_i\}=\mathbb{P}\{X(0)=i\} \quad \textrm{and} \quad  \mathbb{E}\{\Phi^{(m)}B_i\}=\mathbb{P}\{X(0)=i, \Phi^{(m)}=1\}.
\end{align*}
Furthermore, if $\mathbb{P}\{X(t)=i\}>0$ for $t\geq 0$ and $i\in\mathbb{S}$, then by Lemma \ref{lem:lem1},   
\begin{align*}
\lim\limits_{N\rightarrow \infty} \widehat{q}_{N,ij}^{(m)} = \lim\limits_{N\rightarrow \infty}
\frac{N^{-1}\sum_{k=1}^N \Phi_k^{(m)} N_{ij}^{(k)}}{ N^{-1}\sum_{k=1}^N \Phi_k^{(m)} Z_{i}^{(k)}}
=\frac{\mathbb{E}\{\Phi^{(m)}N_{ij}\}}{\mathbb{E}\{\Phi^{(m)} Z_i\}}=q_{ij}^{(m)},
\end{align*}
which in turn implies following (\ref{eq:matq}) that $\lim\limits_{N\rightarrow \infty} \widehat{q}_{N,i}^{(m)}=q_i^{(m)}.$ Moreover, we have
\begin{align*}
\lim\limits_{N\rightarrow \infty} \widehat{s}_{N,i}^{(m)}=\lim\limits_{N\rightarrow \infty}\frac{N^{-1}\sum_{k=1}^N \Phi_k^{(m)} B_i^{(k)}}{N^{-1}\sum_{k=1}^N B_i^{(k)}}=\frac{\mathbb{E}\{\Phi^{(m)} B_i\}}{\mathbb{E}\{B_i\}}=s_i^{(m)},
\end{align*}
where the last equality is due to applying the Bayes' formula for conditional probability, i.e., $\mathbb{P}\{\Phi^{(m)}=1\vert X(0)=i\}=\mathbb{P}\{X(0)=i, \Phi^{(m)}=1\}/\mathbb{P}\{X(0)=i\}.$

\medskip

Finally, we have $\lim\limits_{N\rightarrow \infty} N^{-1}\sum_{k=1}^N B_i^{(k)}=\mathbb{E}\{B_i\}=\pi_i$. All limits hold with probability one, due to the law of large numbers. These justify the claim. \exit 
\end{proof}

%\subsubsection{Asymptotic normality}

\section{Estimation with incomplete information}

%\subsection{Kullback-Leibler divergence}

\subsection{The EM algorithm}
Note that under complete observations, the estimators $\widehat{q}_{N,ij}^{(m)}$, $\widehat{q}_{N,i}^{(m)}$, and $\widehat{s}_{N,i}^{(m)}$ of the distribution parameters of the process $X$ (\ref{eq:mixture}) are given based on knowing the random variable $\Phi_k^{(m)}$ which provides information about the underlying process that drives the $k$th realization of $X$. Under incomplete information, where we only know the realizations $\{\mathbf{X}_k\}$ of $X$, we need to replace $\Phi_k^{(m)}$ by the corresponding estimator $\widehat{\Phi}_k^{(m)}$ based on the observed sample $\{\mathbf{X}_k\}$ of the process.

For this purpose, we apply the EM algorithm, see Dempster et al. \cite{Dempster} and McLachlan and Krishnan \cite{McLachlan} for more details. The first step of the iteration, the $\mathrm{E}-$step, consists of calculating the conditional expectation of the sufficient statistics $\Phi_k^{(m)}N_{ij}^{(k)}$, $\Phi_k^{(m)}N_{i}^{(k)}$, $\Phi_k^{(m)}Z_{i}^{(k)}$, and $\Phi_k^{(m)}B_{i}^{(k)}$, given the sample $\{\mathbf{X}_k\}$. Then, in the $\mathrm{M}-$step, the log-likelihood (\ref{eq:log-like}) is maximized, using the conditional expectation of the sufficient statistics as its observed value. The new estimates of the parameters are given by replacing the statistics in the estimators (\ref{eq:est1})-(\ref{eq:est4}) by their corresponding conditional expectations evaluated at the $\mathrm{E}-$step.

\medskip

The EM algorithm is given below, which follows by an adaptation of \cite{Frydman2005}.

{\setlist[itemize]{leftmargin=*}
\begin{itemize}[label={}]

\item \textbf{Step 1}. Choose initial values of the distribution parameters $\pi_{0,i}$, $q_{0,ij}^{(m)}$, $q_{0,i}^{(m)}$ and $s_{0,i}^{(m)}$ for $i,j=1,\dots,p$, and $m=1,...,M$, all denoted by a vector $\boldsymbol{\theta}_0$. 

\begin{Rem}
Note that the estimator $\widehat{\pi}_i$ of the distribution $\pi_i$ does not get updated at each iteration. It is estimated separately by $\widehat{\pi}_{N,i}$ (\ref{eq:est1}).
\end{Rem}

\item \textbf{Step 2 (E-step)} For the $k$th realization $\mathbf{X}_k$ of $X$, recall that
\begin{align*}
f_{\boldsymbol{\theta}_0}(\Phi_k^{(m)}=1,\mathbf{X}_k)=&f_{\boldsymbol{\theta}_0}(X_k(0)=i)f_{\boldsymbol{\theta}_0}(\Phi_k^{(m)}=1\vert X_k(0)=i)\\
& \times f_{\boldsymbol{\theta}_0}(\mathbf{X}_k\vert \Phi_k^{(m)}=1,X_k(0)=i)\\
=& \big(\pi_{0,i} s_{0,i}^{(m)}\big)^{B_i^{(k)}} f_{\boldsymbol{\theta}_0}(\mathbf{X}_k\vert \Phi_k^{(m)}=1,X_k(0)=i)\\
=& \big(\pi_{0,i} s_{0,i}^{(m)}\big)^{B_i^{(k)}}\prod_{i=1}^p \prod_{j\neq i}^p \big(q_{0,ij}^{(m)}\big)^{N_{ij}^{(k)}} e^{-q_{0,ij}^{(m)} Z_i^{(k)}}.
\end{align*}
For $1\leq m\leq M$, compute the probability $\widehat{\Phi}_k^{(m)}=\mathbb{E}_{\boldsymbol{\theta}_0}\{\Phi_k^{(m)} \big \vert \mathbf{X}_k\}$ that $X$ comes from regime $m$.
Given that $\Phi_k^{(m)}$ is a Bernoulli random variable, 
\begin{align*}
\widehat{\Phi}_k^{(m)}=&f_{\boldsymbol{\theta}_0}(\Phi_k^{(m)}=1\vert \mathbf{X}_k) 
=\frac{f_{\boldsymbol{\theta}_0}(\Phi_k^{(m)}=1, \mathbf{X}_k)}{\sum_{m=1}^M f_{\boldsymbol{\theta}_0}(\Phi_k^{(m)}=1, \mathbf{X}_k)}\\[6pt]
=&\frac{ \prod_{i=1}^p (\pi_{0,i} s_{0,i}^{(m)})^{B_i^{(k)}}\prod_{i=1}^p \prod_{j\neq i}^p \big(q_{0,ij}^{(m)}\big)^{N_{ij}^{(k)}} e^{-q_{0,ij}^{(m)} Z_i^{(k)}} }{  \sum_{m=1}^M \prod_{i=1}^p (\pi_{0,i} s_{0,i}^{(m)})^{B_i^{(k)}}\prod_{i=1}^p \prod_{j\neq i}^p \big(q_{0,ij}^{(m)}\big)^{N_{ij}^{(k)}} e^{-q_{0,ij}^{(m)} Z_i^{(k)}}  }.
\end{align*}
\begin{Rem}
It is straightforward to check that the probability $\widehat{\Phi}_k^{(m)}$ satisfies the constraint $\sum_{m=1}^M \widehat{\Phi}_k^{(m)}=1$, improving the result given in \cite{Frydman2005}. 
\end{Rem}

Then, for $i,j=1,\dots,p$ and $m=1,\dots,M$, compute the conditional expectation of the sufficient statistics $\Phi_k^{(m)}N_{ij}^{(k)}$, $\Phi_k^{(m)}N_{i}^{(k)}$, $\Phi_k^{(m)}Z_{i}^{(k)}$, and $\Phi_k^{(m)}B_{i}^{(k)}$, given the sample $\{\mathbf{X}_k\}$. Note that the random variables $N_{ij}^{(k)}$, $N_{i}^{(k)}$, and $B_{i}^{(k)}$ are all adapted to the information set generated by $\{\mathbf{X}_k\}$.
\begin{align*}
\mathbb{E}_{\boldsymbol{\theta}_0}\big\{\Phi_k^{(m)}N_{ij}^{(k)} \big\vert \mathbf{X}_k\big\}=&\widehat{\Phi}_k^{(m)}N_{ij}^{(k)}\\[6pt]
\mathbb{E}_{\boldsymbol{\theta}_0}\big\{\Phi_k^{(m)}N_{i}^{(k)} \big\vert \mathbf{X}_k\big\}=&\widehat{\Phi}_k^{(m)}N_{i}^{(k)}\\[6pt]
\mathbb{E}_{\boldsymbol{\theta}_0}\big\{\Phi_k^{(m)}Z_{i}^{(k)} \big\vert \mathbf{X}_k\big\}=&\widehat{\Phi}_k^{(m)}Z_{i}^{(k)}\\[6pt]
\mathbb{E}_{\boldsymbol{\theta}_0}\big\{\Phi_k^{(m)}B_{i}^{(k)} \big\vert \mathbf{X}_k\big\}=&\widehat{\Phi}_k^{(m)}B_{i}^{(k)}.
\end{align*}

\item \textbf{Step 3 (M-step)} Compute the new values $\pi_{1,i}$, $q_{1,ij}^{(m)}$, $q_{1,i}^{(m)}$ and $s_{1,i}^{(m)}$ for $i,j=1,\dots,p$, and $m=1,...,M$, using (\ref{eq:est2}), (\ref{eq:est3}), and (\ref{eq:est4}) by
\begin{eqnarray*}
q_{1,ij}^{(m)}&=&\;\frac{\sum_{k=1}^N \widehat{\Phi}_k^{(m)} N_{ij}^{(k)}}{ \sum_{k=1}^N \widehat{\Phi}_k^{(m)} Z_{i}^{(k)}},\label{eq:est2b}\\[6pt]
q_{1,i}^{(m)} &=&\frac{\sum_{k=1}^N \widehat{\Phi}_k^{(m)} N_i^{(k)}}{\sum_{k=1}^N \widehat{\Phi}_k^{(m)} Z_{i}^{(k)}}, \label{eq:est3b} \\[6pt]
s_{1,i}^{(m)} &=&\frac{\sum_{k=1}^N \widehat{\Phi}_k^{(m)} B_i^{(k)}}{\sum_{k=1}^N  B_i^{(k)}}. \label{eq:est4b}
\end{eqnarray*}
Notice that we have replaced the random variable $\Phi_k^{(m)}$ by its respective estimate $\widehat{\Phi}_k^{(m)}$. Stack all the updated estimates into a new vector $\boldsymbol{\theta}_1$.

\begin{Rem}
For the restricted mixture, the new update $s_{1,i}^{(m)}$ is the same as above. However, the updates $\psi_{1,i}^{(m)}$ and $q_{1,ij}$ are given for $j\neq i$ by
\begin{align*}
q_{1,i}=&\frac{\sum_{k=1}^N \widehat{\Phi}_k^{(M)} N_i^{(k)}}{\sum_{k=1}^N \widehat{\Phi}_k^{(M)} Z_i^{(k)}},\\[6pt]
\psi_{1,i}^{(m)}=& \frac{\sum_{k=1}^N \widehat{\Phi}_k^{(m)} N_i^{(k)}}{q_{1,i} \sum_{k=1}^N \widehat{\Phi}_k^{(m)} Z_i^{(k)}}, \\[6pt]
\pi_{1,ij}=&\frac{\sum_{k=1}^N N_{ij}^{(k)}}{\sum_{k=1}^N N_i^{(k)}}, \nonumber\\[6pt]
q_{1,ij}=&\pi_{1,ij} q_{1,i},  \nonumber\\[6pt]
q_{1,ij}^{(m)}=&\psi_{1,i}^{(m)}q_{1,ij}.
\end{align*}

\end{Rem}

\item \textbf{Step 4} Stop if the convergence criterion is achieved. Otherwise, return to Step 1 by replacing $q_{0,ij}^{(m)}$, $q_{0,i}^{(m)}$ and $s_{0,i}^{(m)}$ for $i,j=1,\dots,p$, and $m=1,...,M$, correspondingly by $q_{1,ij}^{(m)}$, $q_{1,i}^{(m)}$ and $s_{1,i}^{(m)}$ for $i,j=1,\dots,p$, and $m=1,...,M$.

For example, the convergence criterion is satisfied when the difference between the updated value $\boldsymbol{\theta}_1$ of each parameter and its previous value $\boldsymbol{\theta}_0$ is less than a specified small positive number, say $\epsilon$, i.e., $\vert\vert \boldsymbol{\theta}_1 - \boldsymbol{\theta}_0 \vert \vert \leq \epsilon$. 

\end{itemize}
}

\section{Simulation study}

To test the performance of the estimation method, we use Monte Carlo simulation to generate sample paths of the mixture process (\ref{eq:mixture}) for a given true values of the distribution parameter. Based on the simulated sample paths, we attempt to estimate the true distribution parameter values using the EM algorithm. 

For this purpose, we assume that the mixture process $X$ (\ref{eq:mixture}) defined on the state space $\mathbb{S}=\{1,2,3\}$ is a mixture of two Markov jump processes $X^{(m)}$, $m=1,2,$ whose intensity matrices $\mathbf{Q}^{(m)}$, $m=1,2,$ can be written as
\begin{align}\label{eq:intensity}
\mathbf{Q}^{(m)}=\textrm{diag}\big(q_1^{(m)},q_2^{(m)},q_3^{(m)}\big)\big(\boldsymbol{\Pi}^{(m)}-\mathbf{I}\big),
\end{align}
where $q_i^{(m)}$ is the exit rate from state $i$ in the $m$'th Markov process $X^{(m)}$, $\textrm{diag}\big(q_1^{(m)},q_2^{(m)},q_3^{(m)}\big)$ is the diagonal matrix, $\boldsymbol{\Pi}^{(m)}$ is the transition matrix of a discrete time Markov chain $Z^{(m)}$ embedded in a continuous Markov process governed by $\mathbf{Q}^{(m)}$, and $\mathbf{I}$ is an identity matrix. Let $\mathbf{q}^{(m)}=\big(q_1^{(m)},q_2^{(m)},q_3^{(m)}\big),$ $m=1,2.$ Expression (\ref{eq:intensity}) suggests the way to carry out the simulation.

\subsection{Specification of the true parameters of the mixture}
Parameter values set for the simulation, which include the initial distribution $\boldsymbol{\pi}$ of starting the process, the exit rates from states $\mathbf{q}_i^{(1)}$ and $\mathbf{q}_i^{(2)}$, $i=1,2,3$, and the speed regime probabilities $\mathbf{s}^{(1)}$ and $\mathbf{s}^{(2)}$ are presented in the following Table. 

\begin{table}[ht]
\centering % used for centering table
\begin{tabular}{c c c c c c} % centered columns (3 columns)
\hline\hline %inserts double horizontal lines
State (i) & $\pi_i$ & $q_i^{(1)}$ & $q_i^{(2)}$ & $s_i^{(1)}$ & $s_i^{(2)}$ \\ [0.5ex] % inserts table
%heading
\hline % inserts single horizontal line
1 &  1/3   &  1/3   & 1/2   & 0.5   & 0.5 \\ % inserting body of the table
2 &  1/3   &   2/5  & 2/5 & 0.25 & 0.75 \\
3 &   1/3  &  1/2   & 1/3   & 0.75 & 0.25 \\[1ex] % [1ex] adds vertical space
\hline %inserts single line
\end{tabular}
\caption{Exit rates and switching probability. } % title of Table
\label{table:switch} % is used to refer this table in the text
\end{table}

The transition matrices of the embedded Markov chains $Z^{(1)}$ and $Z^{(2)}$ are given respectively by

\begin{equation*}
\boldsymbol{\Pi}^{(1)} =
 \left(\begin{array}{cccc}
0 & 0.6 & 0.4 \\
 0.5 & 0 & 0.5\\
 0.4 & 0.6 & 0
\end{array}\right) \quad \textrm{and} \quad 
\boldsymbol{\Pi}^{(2)} =
 \left(\begin{array}{cccc}
 0 & 0.8 & 0.2 \\
 0.5 & 0 & 0.5\\
 0.2 & 0.8 & 0
 \end{array}\right).
\end{equation*}

\subsection{Simulation of the mixture sample paths on $[0,T]$}

From the Monte Carlo method discussed in Section 2, the sample paths of the mixture process can be generated using the following steps.

{\setlist[itemize]{leftmargin=*}
\begin{itemize}[label={}]

\item \textbf{Step 1} Draw at random an initial state $X_0=i_0$ with the distribution $\boldsymbol{\pi}$ on the states $1,2,3$ using the construction (\ref{eq:initialcond}). 

\item \textbf{Step 2} Given the initial state $i_0$, draw using the construction (\ref{eq:initialregime}) the regime indicator from the Bernoulli distribution with the success probability equal to $s_{i_0}$, where success corresponds to regime $\mathbf{Q}^{(1)}$.

\item \textbf{Step 3}\footnote{For this step, we refer among others to Sigman \cite{Sigman}.} Given initial state $i_0$ and regime $m$, that is $Z_0^{(m)}=i_0$, simulate using the recursive equation (\ref{eq:funcF}) and (\ref{eq:chain}) $Z_1^{(m)}$ as follows. Draw $V_1\sim U(0,1)$.

{\setlist[itemize]{leftmargin=*}
\begin{itemize}[label={}]

\item If $Z_0^{(m)}=1$, use the first row of $\boldsymbol{\Pi}^{(m)}$
\begin{align}
\textrm{If} \; V_1 \; \leq& \; \pi_{12}^{(m)},\; \textrm{set} \; Z_1^{(m)}=2 \nonumber \\
\textrm{If} \; V_1 \; >& \; \pi_{12}^{(m)},\; \textrm{set} \; Z_1^{(m)}=3. \label{eq:step1}
\end{align} 

\item If $Z_0^{(m)}=2$, use the second row of $\boldsymbol{\Pi}^{(m)}$
\begin{align}
\textrm{If} \; V_1 \; \leq& \; \pi_{21}^{(m)},\; \textrm{set} \; Z_1^{(m)}=1 \nonumber \\
\textrm{If} \; V_1 \; >& \; \pi_{21}^{(m)},\; \textrm{set} \; Z_1^{(m)}=3. 
\end{align} 

\item If $Z_0^{(m)}=3$, use the third row of $\boldsymbol{\Pi}^{(m)}$
\begin{align}
\textrm{If} \; V_1 \; \leq& \; \pi_{31}^{(m)},\; \textrm{set} \; Z_1^{(m)}=1 \nonumber \\
\textrm{If} \; V_1 \; >& \; \pi_{31}^{(m)},\; \textrm{set} \; Z_1^{(m)}=2.  \label{eq:step3}
\end{align} 

More generally, given $Z_{j-1}^{(m)}=i_j$, $1\leq j\leq J$, to simulate the value of $Z_j^{(m)}$ draw $V_j$ from $U(0,1)$ independently of $V_1,V_2,\dots,V_{j-1}$ and use (\ref{eq:step1})-(\ref{eq:step3}) with $Z_0^{(m)}$ replaced by $Z_{j-1}^{(m)}$, $V_1$ by $V_j$, and $Z_1^{(m)}$ by $Z_{j-1}^{(m)}$. Repeating this procedure $J$ times will generate a sample path $\{Z_0^{(m)}=i_0, Z_1^{(m)}=i_1,\dots,Z_J^{(m)}=i_J\}$. 

\end{itemize} 
}

\item \textbf{Step 4} Simulate the waiting times indicated by the path obtained in Step 3. Since waiting times in states are independent of each other and have exponential distributions with state dependent parameters we simulate them by using independent draws from the exponential distributions corresponding to the sequence of states in the simulated path. Denote the waiting time in state $i_j$ by $S_{i_j}^{(m)}$. Then $S_{i_j}^{(m)}$ has exponential distribution with parameter $q_{i_j}^{(m)}$, i.e., $S_{i_j}^{(m)}=-\log W_{i_j}/q_{i_j}^{(m)}$ with independently drawn $W_{i_j}\sim U(0,1)$. We sequentially generate the draws from the exponential distributions until the first time their sum exceed time $T$, that is when the epoch time, see (\ref{eq:epoch2}), $T_j:=\sum_{k=0}^{j-1}S_{i_k}^{(m)}<T$ and $T_{j+1}>T$. Combining information from Step 3 with the present one gives a sample path of $X^{(m)}$, see the recursive equations (\ref{eq:funcF})-(\ref{eq:simulated}) for details. This sample path is of the form $(Z_0^{(m)}=i_0,S_{i_0}^{(m)}, Z_1^{(m)}=i_1,S_{i_1}^{(m)},\dots, Z_j^{(m)}=i_J,S_{i_J,c}^{(m)})$, where $i_J$ is the last observed state before $T$ and $S_{i_J,c}^{(m)}$ is the censored duration in state $i_J$ by $T$. 

\item \textbf{Step 5} Stop if $N=\textrm{the number of realizations from the mixture process},$ which set to be equal to $20,000$. Otherwise go back to Step 1.

The $20,000$ realizations of the mixture process can be used as an input to the EM algorithm for estimation of the distribution parameters.

\end{itemize} 
}

Figure \ref{fig:Fig1} displays five randomly sampled sample paths of the mixture process. We see that on a given observation time $[0,t]$, $t>0$, each state contains a mixture of two Markov jump processes $X^{(1)}$ and $X^{(2)}$ moving at different speed.    

\subsection{The EM estimation results}

We generate $N=20,000$ independent sample paths of the mixture process $X$ making transitions on the interval $[0,T]$, with $T=100$. Simulation results on five randomly selected sample paths are displayed in Figure \ref{fig:Fig1}. The initial parameter values for $\boldsymbol{\pi}_0$ is set to be equal to the estimate $\widehat{\boldsymbol{\pi}}$, see Table 2 below, while the switching probabilities $\mathbf{s}_0^{(1)}$ and $\mathbf{s}_0^{(2)}$ are chosen randomly on $[0,1]$, whereas the transition matrices $\boldsymbol{\Pi}^{(1)}$ and $\boldsymbol{\Pi}^{(2)}$ are set to be equal to the transition matrix 

\begin{equation}\label{eq:initialmatrx}
\boldsymbol{\Pi} =
 \left(\begin{array}{cccc}
0  & 0.6904  &  0.3096 \\
 05032  &0  &  0.4968 \\
 0.3067  &  0.6933  &  0
\end{array}\right),
\end{equation}
assuming that $\{\mathbf{X}_k\}$ was generated by a simple Markov process with exit rate $\mathbf{q}=(0.3902,0.3902,0.4)$. Based on the sample paths, we obtain:

\begin{table}[ht]
\centering % used for centering table
\begin{tabular}{c c c c c c} % centered columns (3 columns)
\hline\hline %inserts double horizontal lines
State (i) & $\widehat{\pi}_i$ & $\widehat{q}_i^{(1)}$ & $\widehat{q}_i^{(2)}$ & $\widehat{s}_i^{(1)}$ & $\widehat{s}_i^{(2)}$ \\ [0.5ex] % inserts table
%heading
\hline % inserts single horizontal line
1 & 0.3352 & 0.3277   & 0.4930   & 0.4913   & 0.5087 \\ % inserting body of the table
2 & 0.3351 & 0.3924   & 0.3939   & 0.2437   & 0.7563 \\
3 & 0.3297 & 0.4978   & 0.3274   & 0.7545   & 0.2455 \\[1ex] % [1ex] adds vertical space
\hline %inserts single line
\end{tabular}
\caption{Estimates of $\pi_i$, $q_i^{(m)}$ and $s_i^{(m)}$, $m=1,2,$ under unrestricted model. } % title of Table
\label{table:estimate2} % is used to refer this table in the text
\end{table}

It is straightforward to check that $\widehat{s}_i^{(1)} + \widehat{s}_i^{(2)}=1$ for all $i=1,\dots,p$, see (\ref{eq:portion}).

\noindent The estimate of the transition matrix $\boldsymbol{\Pi}^{(1)}$ of the Markov chain  $Z^{(1)}$ is given by
\begin{equation*}
\widehat{\boldsymbol{\Pi}}^{(1)} =
 \left(\begin{array}{ccc}
0 & 0.5927 & 0.4073 \\
 0.5082 & 0 & 0.4918\\
 0.3982 & 0.6018 & 0
\end{array}\right),
\end{equation*}
whereas the estimate of the transition matrix $\boldsymbol{\Pi}^{(2)}$ of $Z^{(2)}$ is found to be
\begin{equation*}
\widehat{\boldsymbol{\Pi}}^{(2)} =
 \left(\begin{array}{ccc}
0 & 0.7974 & 0.2026 \\
 0.4992 & 0 & 0.5008\\
 0.2032 & 0.7968 & 0
\end{array}\right).
\end{equation*}

From the EM estimation outcomes, we observe that the estimates are reasonably close enough to the true values of the distribution parameters.

\subsubsection{Estimation based on the restricted mixture}

The EM estimation of the distribution parameters is based on the restricted model $\mathbf{Q}^{(1)}=\boldsymbol{\Psi}\mathbf{Q}^{(2)}$, with $\boldsymbol{\Psi}=\textrm{diag}(\psi_1, \dots,\psi_p)$. The initial condition for $\mathbf{Q}^{(2)}$ is defined by the matrix $\mathbf{Q}$ (\ref{eq:initialmatrx}), whilst $\psi_i$, the initial distribution $\pi_i$, $s_i^{(1)}$ and $s_i^{(2)}$, for $i=1,\dots,p$, are all chosen randomly on the unit interval $[0,1]$.

\begin{table}[ht]
\centering % used for centering table
\begin{tabular}{c c c c c c} % centered columns (3 columns)
\hline\hline %inserts double horizontal lines
State (i) & $\widehat{\pi}_i$ & $\widehat{q}_i^{(1)}$ & $\widehat{q}_i^{(2)}$ & $\widehat{s}_i^{(1)}$ & $\widehat{s}_i^{(2)}$ \\ [0.5ex] % inserts table
%heading
\hline % inserts single horizontal line
1 & 0.3352  & 0.3302   & 0.4977   & 0.5446   & 0.4554 \\ % inserting body of the table
2 & 0.3351 &  0.3909   & 0.3956   & 0.2690   & 0.7310 \\
3 & 0.3297 &  0.4922   & 0.3256   & 0.7718   & 0.22282 \\[1ex] % [1ex] adds vertical space
\hline %inserts single line
\end{tabular}
\caption{Estimates of $\pi_i$, $q_i^{(m)}$ and $s_i^{(m)}$, $m=1,2,$ under restricted model. } % title of Table
\label{table:estimate2b} % is used to refer this table in the text
\end{table}

It is straightforward to check that $\widehat{s}_i^{(1)} + \widehat{s}_i^{(2)}=1$ for all $i=1,\dots,p$, see (\ref{eq:portion}).

\noindent The estimate of the transition matrix $\boldsymbol{\Pi}^{(1)}$ of the Markov chain $Z^{(1)}$ is
\begin{equation*}
\widehat{\boldsymbol{\Pi}}^{(1)} =
 \left(\begin{array}{ccc}
0 & 0.6904 & 0.3096 \\
0.5032 & 0 & 0.4968\\
 0.3067 & 0.6933 & 0
\end{array}\right),
\end{equation*}
whereas the estimate of the transition matrix $\boldsymbol{\Pi}^{(2)}$ of $Z^{(2)}$ is found to be
\begin{equation*}
\widehat{\boldsymbol{\Pi}}^{(2)} =
 \left(\begin{array}{ccc}
0 & 0.6904 & 0.3096 \\
 0.5032 & 0 & 0.4968\\
 0.3067 & 0.6933 & 0
\end{array}\right),
\end{equation*}
whilst the estimate of the speed reference variable is given by $\widehat{\boldsymbol{\Psi}}=\textrm{diag}(0.6635, 0.9881, 1.5118)$ satisfying the constraint $\widehat{\mathbf{Q}}^{(1)}=\widehat{\boldsymbol{\Psi}}\widehat{\mathbf{Q}}^{(2)}.$
Notice that the EM estimation $\widehat{\boldsymbol{\Pi}}^{(m)}$, $m=1,2,$ is equal to the transition matrix $\boldsymbol{\Pi}$ (\ref{eq:initialmatrx}) of the Markov chain, see Remark \ref{rem:restestm}. As we can see, the EM estimations for the unrestricted mixture outperform that of for the restricted mixture model \cite{Frydman2005}.  

In the section below a statistical test is performed to compare the statistical significance of the Markov model against the Markov mixture model, and the restricted mixture model \cite{Frydman2005} against the unrestricted model at a certain level.

\subsection{Likelihood ratio test}

To test the hypothesis $H_0: q_{i,j}^{(m)}=q_{ij}$ for $i,j=1,\dots,p$ and $m=1,\dots,M$ that the simulated processes is driven by a Markov jump process against the alternative hypothesis $H_1: q_{i,j}^{(m)}\neq q_{ij}$ that it is a mixture of $M$ Markov jump processes, we apply the likelihood ratio test by adapting the one presented in \cite{Frydman2005}. The test statistic is described as follows. Under the $H_0$ hypothesis, the likelihood of observing $N-$independent realizations of the sample paths of $X$ is given by
\begin{align*}
L_{\textrm{Markov}}(\boldsymbol{\pi},\mathbf{Q})=\prod_{k=1}^N f_{\boldsymbol{\theta}_0}(\mathbf{X}_k)= \prod_{k=1}^N\prod_{i=1}^p (\pi_i)^{B_i^{(k)}} \prod_{i=1}^p \prod_{j\neq i}^p (q_{ij})^{N_{ij}^{(k)}} e^{-q_{ij} Z_i^{(k)}},
\end{align*}
where $\mathbf{Q}$ denotes the intensity matrix of a Markov jump process and $Z_i^{(k)}$ is the total time the $k-$th realization $X_k$ stays in state $i$. For the mixture model, the likelihood contribution of observing the sample paths $\mathbf{X}_k$ is given by $f_{\boldsymbol{\theta}_1}(\mathbf{X}_k)=\sum_{m=1}^M f_{\boldsymbol{\theta}_1}(\mathbf{X}_k, \Phi_k^{(m)}) $. Thus, the likelihood under the mixture model is given by 
\begin{align}\label{eq:likelihood2}
&L_{\textrm{Mixture}}\big(\boldsymbol{\pi},\mathbf{Q}^{(m)}, \mathbf{S}^{(m)}, m=1,\dots,M \big)=\prod_{k=1}^N f_{\boldsymbol{\theta}_1}(\mathbf{X}_k)  \nonumber\\ 
&\hspace{2cm}=\prod_{k=1}^N \Big(\sum_{m=1}^M \prod_{i=1}^p (s_i^{(m)}\pi_i)^{B_i^{(k)}} \prod_{i=1}^p \prod_{j\neq i}^p (q_{ij}^{(m)})^{N_{ij}^{(k)}} e^{-q_{ij}^{(m)} Z_i^{(k)}}\Big),
\end{align}
where the product is for all realizations. The likelihood ratio statistic is given by
\begin{align}\label{eq:LRT}
\boldsymbol{\Lambda}_1=\frac{L_{\textrm{Markov}}(\widehat{\boldsymbol{\pi}},\widehat{\mathbf{Q}})}{L_{\textrm{Mixture}}\big(\widehat{\boldsymbol{\pi}},\widehat{\mathbf{Q}}^{(m)}, \widehat{\mathbf{S}}^{(m)}, m=1,\dots,M \big)},
\end{align}
where $\widehat{\boldsymbol{\pi}}$,  $\widehat{\mathbf{Q}}$, $\widehat{\mathbf{Q}}^{(m)}$ and $\widehat{\mathbf{S}}^{(m)}$, with $m=1,\dots,M$, are the MLEs of $\boldsymbol{\pi}$,  $\mathbf{Q}$, $\mathbf{Q}^{(m)}$ and $\mathbf{S}^{(m)}$, with $m=1,\dots,M$. Notice that as the two likelihoods share the same $\boldsymbol{\pi}$, the term $\prod_{k=1}^N\prod_{i=1}^p (\pi_i)^{B_i^{(k)}}$ gets canceled out from the likelihood ratio statistic. The statistic (\ref{eq:LRT}) is calculated for $\mathbf{Q}^{(m)}=\boldsymbol{\Psi}\mathbf{Q}$ and for the unrestricted model. 

To test the significance of the restricted mixture model against the unrestricted mixture, we consider the following likelihood ratio test statistic
\begin{align}\label{eq:LRT2}
\boldsymbol{\Lambda}_2=\frac{L_{\textrm{RestMixture}}\big(\widehat{\boldsymbol{\pi}},\widehat{\mathbf{Q}}, \boldsymbol{\Psi}^{(m)}, \widehat{\mathbf{S}}^{(m)}, m=1,\dots,M \big)}{L_{\textrm{Mixture}}\big(\widehat{\boldsymbol{\pi}},\widehat{\mathbf{Q}}^{(m)}, \widehat{\mathbf{S}}^{(m)}, m=1,\dots,M \big)}.
\end{align}

Given that the entries of intensity matrices $\mathbf{Q}$, $\mathbf{Q}^{(m)}$ and $\mathbf{S}^{(m)}$ respectively satisfy the constraint (\ref{eq:matq}) and (\ref{eq:portion}), it is known by standard theory that under the null hypothesis $-2\ln \boldsymbol{\Lambda}_1$ has $\chi_d^2-$distribution with d.f. $d=p^2(m-1)$. Similarly, with the same arguments, under the null hypothesis $H_0:\psi_i^{(m)}=1,1\leq i\leq p, 1\leq m\leq M-1$ versus $H_1: \textrm{at least one}\;\psi_i^{(m)}\neq 1$, $1\leq i\leq p, 1\leq m\leq M-1$, the test statistic $-2\ln \boldsymbol{\Lambda}_2$ has $\chi_d^2-$distribution with d.f. $d=p(p-1)(m-1)$.

On account that (\ref{eq:likelihood2}) admits no closed-form solution for the MLE estimates of $\boldsymbol{\pi}$,  $\mathbf{Q}^{(m)}$ and $\mathbf{S}^{(m)}$, with $m=1,\dots,M$, we therefore necessarily use the EM estimates in the test statistic (\ref{eq:LRT}) given that the estimation results are reasonably close enough to the actual parameter values of the mixture distribution.

Based on the MLE estimations, the likelihood ratio statistic $-2\ln \boldsymbol{\Lambda}_1$ for comparing the Markov model against the alternative restricted mixture model and unrestricted mixture model has in each case the observe value $2.4311+04$ and $8.0193e+03$, respectively. On the other hand, we have  $-2\ln \boldsymbol{\Lambda}_2=1.6291e+04$. Each alternative is found to be significant at the level $\alpha=5\%$. We therefore conclude that unrestricted mixture model is found to be statistically more significant at the level $\alpha=5\%$ than the Markov and restricted mixture models.  

\section{Conclusions}

We have developed tractable construction of a continuous-time stochastic process based on a finite mixture of right-continuous Markov jump processes moving at different speeds on the same finite state space. As discussed in more details in Frydman and Schuermann \cite{Frydman2008} and Surya (\cite{Surya2018b}, \cite{Surya2018}), that unlike the underlying Markov processes the mixture itself lacks stationarity and the Markov property. 

Monte Carlo method for simulating the process was discussed along with proving distributional equivalence between the simulated process and the theoretical corresponding process. Maximum likelihood estimation was presented for complete and incomplete information. Under complete information, consistent estimators of the distribution parameters were obtained in closed form in terms of sufficient statistics of the process. The EM estimation was proposed for incomplete information knowing only the sample paths of the process. 

Based on Monte Carlo simulation, the EM estimations for the unrestricted mixture were shown to be close enough to the actual value of the distribution parameters, and is found to be statistically significant based on the likelihood ratio test statistic at the level $\alpha=5\%$ compared to the Markov model and restricted  mixture model \cite{Frydman2005}. The results presented in this paper offer appealing features for various applications, for instance in estimating the distribution of first exit time to absorbing state of the mixture process, see for e.g. Surya \cite{Surya2018c}.

\section{Acknowledgments}

Part of the research findings of this paper were presented at the Department of Statistics Seminar of Auckland University on 21 November 2018. The author would like to thank seminar participants and to Simon Harris for the invitation and valuable feedbacks on the work. The author also acknowledges financial support provided by the School of Mathematics and Statistics of Victoria University of Wellington through the school strategic research grant \#220859.

\end{document}